\documentclass[12pt,final]{amsart}
\usepackage{style}
\usepackage{bm}
\usepackage{tikz,mfabacus}

\newcommand{\ak}[2]{\mathcal{AK}_{#2}^{#1}}

\newcommand{\Ei}{\mathcal{E}_{i}}
\newcommand{\Fi}{\mathcal{F}_{i}}

\newcommand{\Cat}[1]{\mathcal{C}_{#1}}

\newcommand{\charge}{\bm\kappa_{\Lambda}}
\newcommand{\la}{\lambda}
\newcommand{\iInd}{\text{$i$}-Ind}
\newcommand{\iRes}{\text{$i$}-Res}

\newcommand{\bQ}{\bm Q}
\newcommand{\std}[2]{\Delta_{#1}^{#2}}
\newcommand{\simple}[2]{L_{#1}^{#2}}

\title[Scopes equivalences for Ariki--Koike algebras as categorical actions]{A note on Scopes equivalences for Ariki--Koike algebras as categorical actions}
\date{}

\author{Alice Dell'Arciprete}
\address[A.D.A.]{Fakult\"at f\"ur Mathematik, Ruhr-Universit\"at Bochum, Universit\"atsstraße 150, D-44780 Bochum.}
\email{alice.dellarciprete@ruhr-uni-bochum.de}

\author{Dinushi Munasinghe}
\address[D.M.]{Department of Mathematics, University of Toronto, 40 St George St, Toronto, ON M5S 2E4}
\email{dinushi.munasinghe@mail.utoronto.ca}
\begin{document}
\maketitle

\begin{abstract}
    A categorical action of a Kac--Moody algebra $\mathfrak{g}$ is built on a category $\mathcal{C}$ decomposed according to the weights $P$ of $\mathfrak{g}$, as well as biadjoint endofunctors $\Ei$ and $\Fi$, abstracting $i$-induction and $i$-restriction, which act on the weight spaces of $\mathcal{C}$ in the same way that the Chevalley generators would act on a regular representation.
Chuang and Rouquier initially developed these notions for $\mathfrak{sl}_2$-categorical actions, using them to prove Broué's abelian defect group conjecture for symmetric groups by establishing derived equivalences between blocks of the same defect. In the setting of general categorical actions Webster later showed that many of these derived equivalences are, in fact, $t$-exact, and that, as a result, such an action can be used to separate weight spaces of a categorical action into a finite number of Morita equivalence classes, where these equivalences also preserve decomposition numbers. The combinatorics of these powerful abstract results were concretely established in the case of Ariki--Koike algebras by the first author in \cite{DA24}, and in this short note we discuss how to translate between the two settings. 
\end{abstract}

\section{Introduction}

Let $\mathfrak{g}$ be an affine Lie algebra with abstract Cartan $\mathfrak{h}$ and simple roots $\alpha_i$ indexed by $i \in I$ (some indexing set). A categorical action of $\mathfrak{g}$ involves a category $\mathcal{C}= \bigoplus_{\nu \in P
} \mathcal{C}_\nu$ decomposed according to the weights $P$ of $\mathfrak{g}$, as well as endofunctors
\begin{align*}
    \mathcal{E}_i &: \mathcal{C}_\nu \rightarrow \mathcal{C}_{\nu+\alpha_i}; \\
     \mathcal{F}_i &: \mathcal{C}_{\nu+\alpha_i} \rightarrow \mathcal{C}_{\nu}.
\end{align*}
These biadjoint functors are an abstraction of induction and restriction functors, and their compositions are required to satisfy $U_q(\mathfrak{g})$-type relations \cite{cautisrigidity}.

Chuang and Rouquier initially developed these notions for $\mathfrak{sl}_2$-categorical actions on abelian categories, using them to prove Broué's abelian defect group conjecture for symmetric groups. That is, they constructed $\mathcal{E}_i$ and $\mathcal{F}_i$ from $i$-induction and $i$-restriction on $\bigoplus_n \Bbbk \mathfrak{S}_n$-modules. By forming a complex of functors categorifying each Weyl group reflection, they established equivalences of derived categories that also preserve decomposition numbers \cite{CR}.

In \cite{rock}, Webster showed, in the setting of general categorical actions, that many of these Chuang--Rouquier equivalences are actually $t$-exact. This establishes that, given such a categorical action, we can separate weight spaces into a finite number of Morita equivalence classes. In this short note we reinterpret the results from the first author's work in terms of categorical actions to see how the decomposition equivalences established combinatorially can be seen in this framework.

The connections we want to emphasize are as follows:

\begin{figure}[h]
  \begin{center}
\scalebox{0.65}{
\begin{tikzpicture}[scale=2]
    \node at (8,9) {$\ak{\Lambda}{}\mmod = \bigoplus_n \ak{\Lambda}{n}\mmod$};
 
    \node at (0,9) {Crystal of Type $\widehat{A}_{e-1}$};
    \node at (0,8) {$\tilde{f}_i$};
    \node at (0,7.5) {$\tilde{e}_i$};
     \node at (3,8) {$\Fi: \mathcal{C}_{\nu+\alpha_i} \rightarrow \mathcal{C}_\nu$};
    \node at (3,7.5) {$\Ei: \mathcal{C}_\nu \rightarrow \mathcal{C}_{\nu + \alpha_i}$};
    \node at (0,8.5) {Crystal Operators};
    \node at (3,9) {Categorical Module $\bigoplus_{\nu \in P} \mathcal{C}_\nu$};
    \node at (3,8.5) {Biadjoint Endofunctors};
    \node at (8,8) {$i-\operatorname{Ind}: \ak{\Lambda}{n}\mmod \rightarrow \ak{\Lambda}{n+1}\mmod$};
    \node at (8,7.5) {$i-\operatorname{Res}:\ak{\Lambda}{n+1}\mmod \rightarrow \ak{\Lambda}{n}\mmod$};
    \draw[->] (8,7) -- node[right]{decategorify} (8,5.5);
    \node at (8,5.25) {$\widehat{\mathfrak{sl}}_e$-module $V(\Lambda)$ of highest weight $\Lambda$};
    \node at (8,5) {Chevalley Generators $\{e_i, f_i\}_{i \in I} \subset \mathfrak{g}$};
    \node at (8,8.5) {$i$-Induce/Restrict};
    \draw[<->] (1,8.5) -- node[above]{\cite{BD}} (1.5,8.5);
      \draw[<->] (5,8.5) -- node[above]{\cite{CR}, \cite{rock}} (6,8.5);
              \node at (0,5.25) {Combinatorial Interpretations for Crystals};
              \node at (0,5) {(Abacus Moves, Young Diagrams)};
              \draw[<->] (0, 5.5) -- (0,7);
\end{tikzpicture}}
\end{center}
    \caption{An Overview}
    \label{fig:placeholder}
\end{figure}
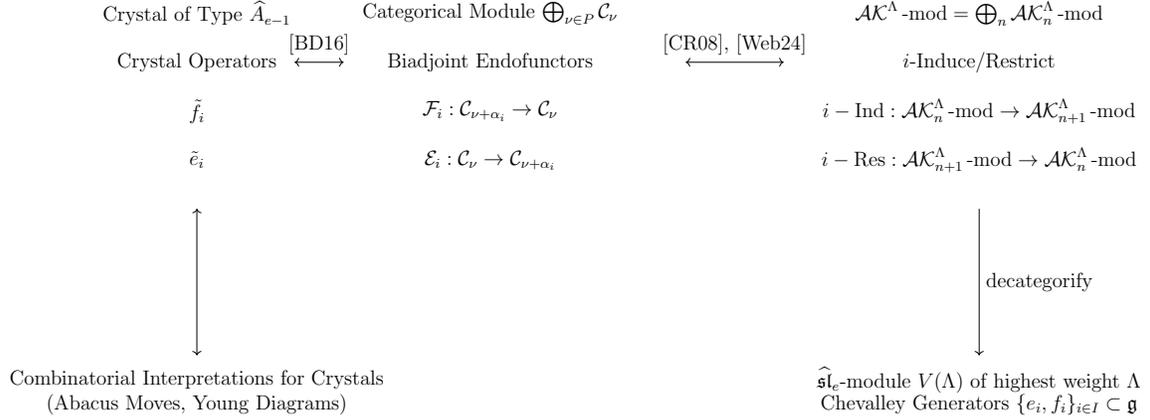
\vspace{10mm}

That is, the powerful theory of categorical actions of Kac--Moody algebras applies directly to the representation theory of Ariki--Koike algebras (which, historically, greatly motivated its development), and so the latter field should benefit from results developed therein.

\subsection*{Acknowledgements} The authors would like to thank Ben Webster for many, many helpful conversations.  The first author is grateful for
financial support from the Alexander von Humboldt Foundation.
The second author would like to thank Maria Chlouveraki, and is grateful for support from the Hellenic Foundation for Research and Innovation (H.F.R.I.) under the Basic Research Financing (Horizontal Support for all Sciences), National Recovery and Resilience Plan (Greece 2.0), Project Number: 15659, Project Acronym: SYMATRAL.

\section{Preliminaries}
Throughout this note, we fix an integer $e$ such that either $e = 0$ or $e\geq 2$. Let $\Gamma$ be the oriented quiver with vertex set $I = \Z$ for $e=0$ or $I = \Z/e\Z$ for $e\geq2$ and with directed edges $i \longrightarrow i + 1$, for all $i \in I$. Thus, $\Gamma$ is the quiver of type $A_{\infty}$ if $e = 0$, and if $e \geq 2$ then it is a cyclic quiver of type $A_e^{(1)}$.

Let $A = (a_{ij})_{i, j \in I}$ be a generalized Cartan matrix associated to $\Gamma$, let $\mathfrak{g}$ be the associated Kac--Moody algebra  generated by $\{ e_i, f_i\}_{i \in I}$, and let $\mathfrak{h}$ be its Cartan subalgebra. Let $\{ \alpha_i\}_{i \in I} \subset \mathfrak{h}^*$ be the (linearly independent) simple roots, $\{ \Lambda_i\}_{i \in I} \subset \mathfrak{h}^*$ be the set of fundamental weights, and let $\langle \cdot , \cdot\rangle$ be the bilinear form determined by $$\langle \alpha_i , \alpha_j \rangle = a_{ij} \qquad \text{and}\qquad\langle \Lambda_i , \alpha_j \rangle = \delta_{i,j}, \quad \text{ for }i, j \in I.$$

Let $P:= \bigoplus_{i \in I} \mathbb{Z}\Lambda_i$ be the weight lattice, generated by the fundamental weights $\Lambda_i$, and let $P_+:=\bigoplus_{i \in I} \mathbb{Z}_{\geq 0}\Lambda_i$ be dominant weight lattice. 

The algebras we are working with are cellular (see \cite{GL} for their definition and properties); given an algebra $\mathcal{A}$ we will denote the standard module indexed by a poset element $\bm\lambda$ as $\std{\mathcal{A}}{\bm\lambda}$, and the simple by $\simple{\mathcal{A}}{\bm\lambda}$. We will refer to the standard modules of the Ariki--Koike algebra as \textbf{Specht modules}, and denote them $\std{\ak{}{}}{\bm\lambda} = S^{\bm\lambda};$ and we will denote the simple modules of the Ariki--Koike algebra by $\simple{\ak{}{}}{\bm\lambda} = D^{\bm\lambda}.$

\section{Representations of Ariki--Koike Algebras}\label{AK}

\subsection*{Ariki--Koike algebras} Let $\mathbb{F}$ be a field. Fix a dominant weight $\Lambda\in P_+$ and
a non-negative integer~$n$.  Set $\ell=\sum_{i\in I}\langle\Lambda,\alpha_i\rangle$. A
\textbf{multicharge} for $\Lambda$ is any sequence of integers
$\charge=(\kappa_1,\dots,\kappa_\ell)\in\Z^\ell$ such that 
$$\langle\Lambda,\alpha_i\rangle=\#\{1\le s\le\ell \ | \ \kappa_s\equiv i\pmod e\}, \text{ for }i\in I,$$ 
where we use the convention that $i\pmod e=i$ if $e=0$. 

Define the
\textbf{quantum characteristic} of $q\in \mathbb{F}$ to be the integer $e$ which is
minimal such that $1+q+\dots+q^{e-1}=0$, and where we set $e=0$ if no such $e$
exists.
Define $\bQ_\Lambda=(Q_{\kappa_1},\dots,Q_{\kappa_\ell})\in \mathbb{F}^{\ell}$, where for an
integer $k\in\Z$ we set
$$Q_k=\begin{cases}q^k,&\text{if }q\ne1,\\
                     k,&\text{if }q=1.
\end{cases}$$
Note that since $\mathbb{F}$ is a field, $\bQ_\Lambda$ depends only on $\Lambda$
and not on the choice of multicharge~$\charge$.

\begin{defn}
Let $\mathbb{F}$ be a field  and $q$ a non-zero element of $\mathbb{F}$. Let $e$ be the quantum characteristic of $q$ and
  $\Lambda\in P_+$ a dominant weight for $\Gamma$. The \textbf{Ariki--Koike algebra} $\ak{\Lambda}{n}:= \ak{}{n}(q, \bQ_\Lambda)$
 
is generated by elements $T_0,T_1 \dotsc, T_{n-1}$, subject to the following relations:

\begin{equation*}
\left.
\begin{array}{crclll}
    &T_0T_1T_0T_1& = &T_1T_0T_1T_0, \\
    &T_iT_{i+1}T_i &= &T_{i+1}T_iT_{i+1} & \text{ for } 1 \leq i \leq n-2,\\
    &T_iT_j &= &T_jT_i & \text{ for } 0 \leq i,j \leq n-1, |i-j|>1, \end{array}\right\} \text{(Braid Relations)}
\end{equation*}
\begin{equation*}
    \begin{array}{crclll}
        & (T_i-q)(T_i+1)&=&0 & \text{ for } 1 \leq i \leq n-1, & \text{(Quadratic Relation)}\\
        &(T_0-Q_{\kappa_1})\dots(T_0-Q_{\kappa_{\ell}})& =&0. & &\text{(Cyclotomic Relation)}
    \end{array} 
\end{equation*}
\end{defn}

\subsection*{Multipartitions}\label{multpar}
A \textbf{partition} of $n$ is defined to be a non-increasing sequence $\lambda = (\lambda_1, \lambda_2, \dots)$ of non-negative integers whose sum is $n$. {The integers $\lambda_b$, for $b\geq1$, are called the \textbf{parts} of $\lambda$. We write $|\lambda| = n$.\\
Since $n < \infty$, there is a $k$ such that $\lambda_b = 0$ for $b > k$ and we write $\lambda = (\lambda_1, \dots , \lambda_k)$. We write $\varnothing$ for the unique empty partition $(0, 0, \dots)$. If a partition has repeated parts, for convenience we group them together with an index. For example, $$(4,3,3,1,1,1,0,0, \dots) = (4,3,3,1,1,1) = (4,3^2,1^3).$$
The \textbf{Young diagram} of a partition $\lambda$ is the subset
$$[\lambda] := \{(b,c) \in \mathbb{N}_{>0} \times \mathbb{N}_{>0
}\text{ } | \text{ } c\leq \lambda_b\}.$$
\begin{defn}\label{multipar}
An $\ell$-\textbf{multipartition} of $n$ is an ordered $r$-tuple $\bm{\lambda} = (\lambda^{(1)}, \dots, \lambda^{(\ell)})$  of partitions such that $$|\bm{\lambda}| := |\lambda^{(1)}|+ \ldots +|\lambda^{(\ell)}| = n.$$
If $\ell$ is understood, we shall just call this a multipartition of $n$.
\end{defn}
We write the unique multipartition of $0$ as $\bm\varnothing$. The \textbf{Young diagram} of a multipartition $\bm\lambda$ is the subset
$$
[\bm{\lambda}]:=\{(b,c,j)\in \mathbb{N}_{>0}\times \mathbb{N}_{>0}\times\{1, \ldots, \ell\} \text{ }|\text{ } c \leq \lambda_b^{(j)}\}.
$$
We may abuse notation by not distinguishing a multipartition from its Young diagram.}
The elements of $[\bm{\lambda}]$ are called \textbf{nodes} of $\bm\lambda$. We say that a node $\mathfrak{n} \in [\bm{\lambda}]$ is \textbf{removable} if $[\bm{\lambda}]\setminus \{\mathfrak{n}\}$ is also the Young diagram of a multipartition. We say that an element $\mathfrak{n} \in \mathbb{N}^2_{>0}\times \{1, \ldots, \ell\}$ is an \textbf{addable node} if $\mathfrak{n} \notin [\bm{\lambda}]$ and $[\bm{\lambda}] \cup \{\mathfrak{n}\}$ is the Young diagram of a multipartition. {Now fix $e \in \{2, 3, \ldots\}\cup \{0\}$. Given an $\ell$-tuple of integers $\bm\kappa=(\kappa_1, \ldots, \kappa_\ell)$, to each node $(b,c,j) \in [\bm\lambda]$ we associate its \textbf{residue} $\mathrm{res}_{\bm\kappa}(b,c,j) = \kappa_j + c-b$ $(\text{mod }e)$ if $e\geq 2$, and $\mathrm{res}_{\bm\kappa}(b,c,j) = \kappa_j + c-b$ if $e=0$. We draw the residue diagram of $\bm\lambda$ by replacing each node in the Young diagram by its residue.}

\begin{example}\label{exsameblock}
Suppose $\ell=3$ and $\bm\kappa=(1,0,2)$. Let $\bm{\lambda}=((1^3),(2,1),(1^2))$ and $\bm\mu=((1),(2,1),(2,1^2))$ be two multipartitions of $8$. If $e=4$, then the $4$-residue diagrams of $[\bm{\lambda}]$ and of $[\bm\mu]$ are
$$
\young(1,0,3) \quad \young(01,3) \quad \young(2,1) \quad \text{ and } \quad \young(1) \quad \young(01,3) \quad \young(23,1,0).
$$
\end{example}
For an $\ell$-multipartition $\bm\lambda$, define the \textbf{residue set} of $\bm\lambda$ to be the multiset $\mathrm{Res}_{\bm\kappa}(\bm\lambda) = \{\mathrm{res}_{\bm\kappa}(\mathfrak{n}) \text{ } | \text{ } \mathfrak{n} \in [\bm\lambda]\}$. 
Notice that in the example above
\begin{equation*}
\mathrm{Res}_{\bm\kappa}(\bm\lambda)= \{0,0,1,1,1,2,3,3\} = \mathrm{Res}_{\bm\kappa}(\bm\mu).
\end{equation*}

\subsection*{$\beta$-numbers and the abacus}\label{betaAbacus}
We may conveniently represent multipartitions on an abacus display. Fix $\charge=(\kappa_1, \ldots, \kappa_\ell)\in \mathbb{Z}^{\ell}$
to be a multicharge of $\ak{\Lambda}{n}$.
\begin{defn}
Let $\bm{\lambda}=(\lambda^{(1)}, \ldots, \lambda^{(\ell)})$ be a multipartition of $n$. For every $i\geq1$ and for every $j \in \{1, \ldots, \ell\}$, we define the \textbf{$\beta$-number} $\beta_i^j$ to be
$$\beta_i^j := \lambda_i^{(j)} + \kappa_j - i.$$
The set $B_{\kappa_j}^j =\{ \beta_1^j, \beta_2^j, \ldots \}$ is the set of $\beta$-numbers (defined using the integer $\kappa_j$) of partition $\lambda^{(j)}$.
\end{defn}

If $e\geq2$, for each set $B_{\kappa_j}^j$, an $e$-abacus for $\lambda^{(j)}$ with respect to $\kappa_j$ is an abacus with $e$ infinite vertical runners which we label from left to right by the elements of $I$. We mark positions on runner $l$ and label them with the integers congruent to $l$ modulo $e$, so that 
position $(x+1)e+l$ lies immediately below position $xe + l$, for each $x$. Now we place a bead at position $\beta_i^j$, for each $i$. Moreover, we say that the bead corresponding to the $\beta$-number $xe+l$ is at \textbf{level} $x$ for $x \in \mathbb{Z}$. Then, we define $\ell^{\charge}_{ij}(\bm{\lambda})$ to be the level of the lowest bead on runner $i$ of the abacus display for $\lambda^{(j)}$ with respect to $\charge$. Note that this bead is the largest element of $B_{\kappa_j}^j$ congruent to $i$ modulo $e$.
If $e=0$, the
abacus has runners and bead positions indexed by the elements of $I=\Z$, so that runner $i\in\Z$ contains either
one bead or no beads.

Hence, we can now define the $e$-\textbf{abacus display}, or the $e$-\textbf{abacus configuration}, for a multipartition $\bm{\lambda}$ with respect to $\charge$ to be the $\ell$-tuple of $e$-abacus displays associated to each component $\lambda^{(j)}$. If it is clear which $e$ we are referring to, we simply say abacus configuration. When we draw abacus configurations we will draw only a finite part of the runners and we will assume that above this point the runners are full of beads and below this point there are no beads. 

\begin{example}
Suppose that $\ell = 3$, $\charge = (3, 1, 1)$ and $\bm{\lambda}= (5,4,1), (3,2^3), (4,3,1))$. Then we have
\begin{align*}
B^1_{3} &= \{7,5,1,-1,-2,-3,-4, \ldots\};\\
B^2_1 &= \{6,4,3,2,-1,-2,-3,-4, \ldots\};\\
B^3_1 &= \{4,2,-1,-3,-4,\ldots\}.
\end{align*}
So, the $4$-abacus display with respect to the multicharge $\charge$ for $\bm{\lambda}$ is

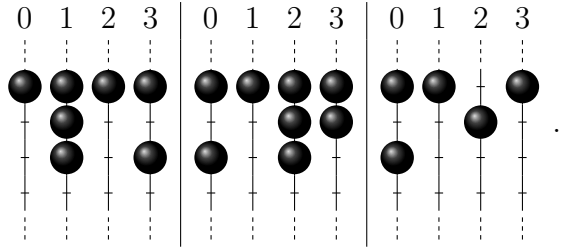
\begin{center}
\begin{tabular}{c|c|c}
$\begin{matrix} 0&1&2&3
\end{matrix}$ & $\begin{matrix} 0&1&2&3
\end{matrix}$ & $\begin{matrix} 0&1&2&3
\end{matrix}$\\
\sabacus(1.5,vvvv,bbbb,nbnn,nbnb,nnnn,vvvv)
&
\sabacus(1.5,vvvv,bbbb,nnbb,bnbn,nnnn,vvvv)
&
\sabacus(1.5,vvvv,bbnb,nnbn,bnnn,nnnn,vvvv)\\
\end{tabular}.
\end{center}
\end{example}

\subsection*{Specht modules and simple modules}\label{specht}
The algebra $\ak{\Lambda}{n}$ is a cellular algebra \cite{djm,GL} with the standard modules indexed by $\ell$-multipartitions of $n$. For each $\ell$-multipartition $\bm\lambda$ of $n$, we define an $\ak{\Lambda}{n}$-module $S^{\bm\lambda}$ called a \textbf{Specht module}; these modules are the standard modules given by the cellular basis of $\ak{\Lambda}{n}$ \cite{djm}. 

\begin{itemize}
    \item When $\ak{\Lambda}{n}$ is semisimple, the Specht modules form a complete set of non-isomorphic irreducible $\ak{\Lambda}{n}$-modules.
    \item 
   
    When $\ak{\Lambda}{n}$ is not semisimple, the Specht modules are not necessarily irreducible. In this case, we need to introduce an important subset of the set of all multipartitions: the subset $\mathcal{K}$ of multipartitions of $n$ which index the simple modules for $\ak{\Lambda}{n}$. The multipartitions in the set $\mathcal{K}$ are called \textit{Kleshchev multipartitions} and they are a subset of $e$-restricted multipartitions (i.e., the set of multipartitions with all components $\lambda^{(j)}$ satisfing  $\lambda_i^{(j)} - \lambda_{i+1}^{(j)} < e$ for every $i\geq1$ and $1 \leq j\leq \ell$). A recursive definition of Kleshchev multipartitions can be found in \cite{Fay07}. In this note, we only need the fact that they index the simple modules of $\ak{\Lambda}{n}$.
    
    Then, for each Kleshchev $\ell$-multipartition $\bm\lambda$, the Specht module $S^{\bm\lambda}$ has an irreducible head $D^{\bm\lambda},$
    and the $D^{\bm\lambda}$ provide a complete set of irreducible modules for $\ak{\Lambda}{n}$ as $\bm\lambda$ ranges over the set of Kleshchev multipartitions of $n$ \cite[Theorem 4.2]{ari01}. 
\end{itemize}

If $\bm\lambda$ and $\bm\mu$ are $\ell$-multipartition of $n$ with $\bm\mu$ Kleshchev, let $d_{\bm\lambda \bm\mu} = [S^{\bm\lambda}:D^{\bm\mu}]$ denote the multiplicity of the simple module $D^{\bm\mu}$ as a composition factor of the Specht module $S^{\bm\lambda}$. The matrix $D = (d_{\bm\lambda \bm\mu})$ is called the \textbf{decomposition matrix} of $\ak{\Lambda}{n}$.

\subsection*{Blocks of Ariki--Koike algebras}

By the cellularity of $\ak{\Lambda}{n}$, each Specht module $S^{\bm\lambda}$ belongs to a unique block of $\ak{\Lambda}{n}$. We therefore abuse notation and say that a multipartition $\bm\lambda$ lies in a block $B$ if the corresponding Specht module $S^{\bm\lambda}$ lies in $B$. Conversely, every block contains at least one Specht module. Hence, classifying the blocks of $\ak{\Lambda}{n}$ amounts to describing the partition of the set of multipartitions. We obtain the following classification of the blocks of $\ak{\Lambda}{n}$.

\begin{thm}\cite[Theorem 2.11]{LM}\label{block_car}
Let $\charge\in\Z^\ell$ be a multicharge of $\ak{\Lambda}{n}$. Suppose $\bm\lambda$ and $\bm\mu$ are $\ell$-multipartitions of $n$. Then, $S^{\bm\lambda}$ and $S^{\bm\mu}$ lie in the same block of $\ak{\Lambda}{n}$ if and only if $\mathrm{Res}_{\charge}(\bm\lambda)=\mathrm{Res}_{\charge}(\bm\mu)$.
\end{thm}

Let $c_i(\bm \lambda)$ denote the number of nodes in $[\bm\la]$ of residue $i\in I$, and let $\alpha (\bm \lambda) := \sum_{i\in I} c_i(\bm \lambda) \alpha_i$. 
\begin{rmk}
    Thus we can index the blocks of the Ariki--Koike algebra by distinct $\alpha$ coming from $\ell$-multipartitions of $n$. If $\alpha = \alpha(\bm \lambda)$, we will say that $B_\alpha$ is the block containing the Specht module indexed by $\bm \lambda.$ In particular, \cref{block_car} can be stated as $\bm\la, \bm\mu$ lie in the same block of $\ak{\Lambda}{n}$ if and only if $\alpha(\bm\lambda)=\alpha(\bm\mu)$.
\end{rmk}

\begin{example}
Continuing Example \ref{exsameblock},  we see that the residue sets of $\bm\lambda$ and $\bm\mu$ are equal. Hence $\bm\lambda$ and $\bm\mu$ lie in the same block of $\ak{\Lambda}{8}$ where $e=4$.   
\end{example}

\begin{defn}\cite[\S 3.1]{Fay07}
    A \textbf{core block}  of $\ak{\Lambda}{n}$ is either any block if $e=0$ or a block $B$ such that every multipartition in the block satisfies the following equivalent conditions:
\begin{enumerate}
\item $\bm{\lambda}$ is a multicore, and there exists a multicharge $\bm{\mathrm{a}} = (a_1, \ldots, a_r)$ such that $a_j \equiv \kappa_j\text{ }(\mathrm{mod}\text{ }e)$ for all $j$ and integers $b_0, \ldots, b_{e-1}$ such that for each $i\in I$ and $j\in \{1, \ldots, \ell\}$, $\ell^{\bm{\mathrm{a}}}_{ij}(\bm{\lambda})$ equals either $b_i$ or $b_i + 1$.
We call such an $e$-tuple $(b_0, \ldots, b_{e-1})$ a \textbf{base tuple} for $\bm{\lambda}$.

\item There is no block of any $\ak{\Lambda}{m}$ with 
the same hub as $B$ and smaller weight than $B$.
\item Every multipartition in $B$ is a multicore.
\end{enumerate}
\end{defn}

\subsection*{Induction and restriction functors}

If $n>1$, then $\ak{\Lambda}{n-1}$ is naturally a subalgebra of $\ak{\Lambda}{n}$, and in fact $\ak{\Lambda}{n}$ is free as an $\ak{\Lambda}{n-1}$-module. So there are well-behaved induction and restriction functors between the module categories of $\ak{\Lambda}{n-1}$ and $\ak{\Lambda}{n}$ given as follows.
For the induction functor we have:
$$\Ind{}{}: \begin{cases} \ak{\Lambda}{n-1}\mmod\longrightarrow\ak{\Lambda}{n}\mmod \\ M\mapsto \Ind{}{}M=M\otimes_{\ak{\Lambda}{n-1}}\ak{\Lambda}{n} \end{cases}$$
For the restriction functor we have:
$$\Res{}{}:\begin{cases} \ak{\Lambda}{n}\mmod\longrightarrow\ak{\Lambda}{n-1} \mmod\\ M\mapsto \Res{}{}M \end{cases}$$

{In \cite{isoAKtoKLR}, it was shown that the Ariki-Koike algebra is isomorphic to a cyclotomic Khovanov-Lauda-Rouquier algebra of type $A$. These are $\Z$-graded algebras and it was further shown in \cite{GradedSpecht} that there is a
corresponding $\Z$-grading on the Specht modules. Using these results one may define graded decomposition numbers $d_{\bm\la\bm\mu}(v)=[S^{\bm\la}:D^{\bm\mu}]_v\in \mathbb{N}[v, v^{-1}]$; we recover the original decomposition numbers by setting $v=1$. For more
details, see \cite{KleshRepThSymHecke}. Given such a grading, it is natural to consider the graded branching rule for the graded Specht modules. In order to state this result, we need to introduce the following notation. We impose a partial order $>$ on the set of nodes of residue $i \in I$ of a multipartition by saying that $(b,c,j)$ is \textbf{above} $(b',c',j')$ (or $(b',c',j')$ is \textbf{below} $(b,c,j)$) if either $j<j'$ or ($j=j'$ and $b<b'$). In this case we write $(b,c,j)>(b',c',j')$.

Let $\bm\la$ be an $\ell$-multipartition of $n$, $i\in I$, $R$ be a removable $i$-node of 
$\bm\la$ and $A$ be an addable $i$-node of $\bm\la$. We set
\begin{equation}\label{EDMUA}
\begin{split}
N_{R}(\bm\la):= &\#\{\text{addable $i$-nodes of $\bm\la$ below $R$}\}
\\
&-\#\{\text{removable $i$-nodes of $\bm\la$ below  $R$}\};
\end{split}
\end{equation}
\begin{equation}\label{EDMUB}
\begin{split}
N^A(\bm\la):=&\#\{\text{addable $i$-nodes of $\bm\la$ above $A$}\}\\
&-\#\{\text{removable $i$-nodes of $\bm\la$ above $A$}\}.
\end{split}
\end{equation} 

\begin{thm}\cite[Main Theorem]{hm11}\cite[Theorem 4.11]{GradedSpecht}\label{gradedbranching}
\begin{itemize}
\item Suppose $\bm \lambda$ is a multipartition of $n-1$, and let $A_1 > \ldots > A_s$ be all the addable nodes of $[\bm \lambda]$. For each $m = 1, \ldots, s$, let $\bm\lambda^{A_m}$ be the multipartition of $n$ with $[\bm\lambda^{A_m}] = [\bm\lambda] \cup \{A_m\}$. Then $\mathrm{Ind}\ S^{\bm\lambda}$ has a filtration
$$\mathrm{Ind}\ S^{\bm\lambda}=V_s \supset \dots \supset V_1 \supset V_0 = \{0\}$$
as a graded $\ak{\Lambda}{n}$-module in which the factors are 
$$ V_m / V_{m-1} = S^{\bm\lambda^{A_m}}\langle N^{{A}_m}(\bm\la)\rangle$$
for all $1 \leq m \leq s$.

\item Suppose $\bm\lambda$ is a multipartition of $n$, and let $R_1 < \ldots < R_t$ be all the removable nodes of $[\bm\lambda]$. For each $m = 1, \ldots, t$, let $\bm\lambda_{R_m}$ be the multipartition of $n-1$ with $[\bm\lambda_{R_m}] = [\bm\lambda] \setminus \{R_m\}$. Then $\mathrm{Res}\ S^{\bm\lambda}$ has a filtration
$$\mathrm{Res}\ S^{\bm\lambda}=V_t \supset \dots \supset V_1 \supset V_0 = \{0\}$$
as a graded $\ak{\Lambda}{n-1}$-module in which the factors are 
$$ V_m / V_{m-1} = S^{\bm\lambda_{R_m}}\langle N_{{R}_m}(\bm\la)\rangle$$
for all $1 \leq m \leq t$.
\end{itemize}
\end{thm}
If a module $M$ has a filtration $M=V_t \supseteq V_{t-1} \supseteq \ldots \supseteq V_{1} \supseteq V_0 = 0$ with $V_m/V_{m-1} \cong M_m\langle d_m\rangle$ for $m = 1, \ldots, t$, and $d_m\in \Z$, we write $$M \sim \sum_{m=1}^{t} M_m \langle d_m\rangle \text{ or } M \sim\sum_{m=1}^{t}v^{d_m} M_m, \qquad \text{for }d_m\in \Z.$$ }

By projecting onto the blocks of $\ak{\Lambda}{n}$ (resp. $\ak{\Lambda}{n-1}$), the induction (resp. restriction) functor decomposes as
$$\Ind{}{} = \bigoplus_{i\in I} \text{\iInd} \qquad \Res{}{} = \bigoplus_{i\in I} \text{\iRes}.$$ 
Then we have the following refinement of \cref{gradedbranching} for $i$-induction (resp. $i$-restriction) functors.

\begin{cor}\cite[Corollary 4.6]{hm11} Let $i\in I$.
\begin{itemize}
\item Suppose $\bm \lambda$ is a multipartition of $n-1$ an, and let $A_1 > \ldots > A_s$ be all the addable $i$-nodes of $[\bm \lambda]$. For each $m = 1, \ldots, s$, let $\bm\lambda^{A_m}$ be the multipartition of $n$ with $[\bm\lambda^{A_m}] = [\bm\lambda] \cup \{A_m\}$. Then $i\text{-}\mathrm{Ind}\ S^{\bm\lambda}$ has a filtration
$$i\text{-}\mathrm{Ind}\ S^{\bm\lambda}=V_s \supset \dots \supset V_1 \supset V_0 = \{0\}$$
as a graded $\ak{\Lambda}{n}$-module in which the factors are 
$$ V_m / V_{m-1} = S^{\bm\lambda^{A_m}}\langle N^{{A}_m}(\bm\la)\rangle$$
for all $1 \leq m \leq s$.

\item Suppose $\bm\lambda$ is a multipartition of $n$, and let $R_1 < \ldots < R_t$ be all the removable $i$-nodes of $[\bm\lambda]$. For each $m = 1, \ldots, t$, let $\bm\lambda_{R_m}$ be the multipartition of $n-1$ with $[\bm\lambda_{R_m}] = [\bm\lambda] \setminus \{R_m\}$. Then $i\text{-}\mathrm{Res}\ S^{\bm\lambda}$ has a filtration
$$i\text{-}\mathrm{Res}\ S^{\bm\lambda}=V_t \supset \dots \supset V_1 \supset V_0 = \{0\}$$
as a graded $\ak{\Lambda}{n-1}$-module in which the factors are 
$$ V_m / V_{m-1} = S^{\bm\lambda_{R_m}}\langle N_{{R}_m}(\bm\la)\rangle$$
for all $1 \leq m \leq t$.
\end{itemize}    
\end{cor}

\subsection*{Weight and hub of a block}

Fix a multicharge $\charge=(\kappa_1, \ldots, \kappa_{\ell})$ of $\ak{\Lambda}{n}$.

\begin{defn}\label{mult_weight}
Let  $\bm{\lambda}= (\lambda^{(1)}, \ldots, \lambda^{(\ell)})$ be a multipartition of $n$. Let 
$\overline{\kappa_j}$ denote $\kappa_j\text{ }(\text{mod }e)$.
Define:
\begin{itemize}
    \item the \textbf{weight} $w(\bm{\lambda})$ of $\bm{\lambda}$ to be the non-negative integer
$$
w(\bm{\lambda}) = \left( \sum\limits_{j=1}^{\ell} c_{\overline{\kappa_j}}(\bm{\lambda})\right) - \frac{1}{2} \sum_{i \in I} (c_{i}(\bm{\lambda})-c_{i+1}(\bm{\lambda}))^2.
$$
\item the \textbf{hub} of $\bm{\lambda}$ to be the collection of integers $(\delta_i(\bm{\lambda}) := \sum_{j=1}^{\ell} \delta_i^j(\bm{\lambda}) \text{ }|\text{ }i \in I)$  where
for each $i\in I$ and $j \in \{1, \ldots, \ell\}$,  
\begin{align*}
\delta_i^j(\bm{\lambda}) = & \#\{\text{removable $i$-nodes of $[\lambda^{(j)}]$}\}\\
 &- \#\{\text{addable $i$-nodes of $[\lambda^{(j)}]$}\},
\end{align*}
\end{itemize}
\end{defn}

Moreover, we want to highlight an important feature of the weight and hub of a multipartition: they are invariants of the block containing $\bm{\lambda}$, and in fact determine the block (see \cite{Fay07} for details). In view of this, we may define the \textbf{weight} (resp. \textbf{hub}) \textbf{of a block} $B$ to be the weight (resp. {hub}) of any multipartition $\bm{\lambda}$ in $B$, and we write $$w(B) = w(\bm{\lambda})\text{ (resp. } \delta_i(B) = \delta_i(\bm{\lambda})).$$

\subsection*{Graded decomposition matrix equivalences}

{We take a moment to recall maps between blocks of Ariki--Koike algebras that are analogous to those defined by Scopes between blocks of symmetric group algebras in \cite{scopesCartanMatrices1991}. These maps were defined by the first author in \cite{DA24}.
Suppose $e\geq 2,$ $i \in I$, and let
$\phi_i \colon \mathbb{Z} \rightarrow \mathbb{Z}$ be the map given by
\begin{equation*}
\phi_i(x)=
\begin{cases}
x + 1 & x \equiv i - 1 \quad(\mathrm{mod}\text{ }e)\\
x - 1 & x \equiv i \quad(\mathrm{mod}\text{ }e)\\
x & \text{otherwise}
\end{cases}.
\end{equation*}

Now suppose $\bm{\lambda}$ is an $\ell$-multipartition, and consider its abacus display with respect to the multicharge $\charge=(\kappa_1, \ldots, \kappa_{\ell})$.
For each $j$, we define a partition $\Phi_i(\lambda^{(j)})$ by replacing each beta-number $\beta$ with $\phi_i(\beta)$.
Equivalently, we simultaneously remove all removable $i$-nodes from $[\lambda^{(j)}]$ and add all addable $i$-nodes of $[\lambda^{(j)}]$. If $i\neq 0$, this is equivalent to swapping runners $i-1$ and $i$ of each component in the abacus display of $\bm\lambda$; {if $i=0$, we swap runners $0$ and $e-1$ and then increase the level of each bead on runner $0$ by $1$ and decrease the level of each bead on runner $e-1$ by $1$.} We define $\Phi_i(\bm{\lambda})$ to be the multipartition $(\Phi_i(\lambda^{(1)}), \ldots,\Phi_i(\lambda^{(\ell)}))$.

In particular, Fayers in \cite[Proposition 4.6]{Fay06} proves that $\Phi_i$ gives a bijection between the set of multipartitions in $B$ of $\ak{\Lambda}{n}$ and the set of multipartitions in $\bar{B}=\Phi_i(B)$ of $\ak{\Lambda}{n-\delta_i(B)}$ having the same weight of $B$.

\begin{thm}\cite[Theorem 5.7]{DA24} \label{AliceMain}
    Fix $i\in I$. Given a block $B$ of $\ak{\Lambda}{n}$ such that $\delta=\delta_i(B)\geq 0$, and $w(B) \leq w(C) + K_i\ell$, where $C$ is the core block of $B$ and 
    \[K_i = \begin{cases}
        b_i-b_{i-1}-1 \quad i \in I \setminus \{ 0\} \\
        b_0-b_{e-1}-2 \quad i=0
    \end{cases}.\]
Suppose that $\bm{\lambda}$ belongs to the block $B$. Set $r = \frac{\delta(\delta-1)}{2}$. Then
\begin{enumerate}
\item[1.] $(i\text{-}\mathrm{Res})^{(\delta)}\ D^{\bm{\lambda}} \sim \sum\limits_{k=0}^{r} |\mathfrak{S}_{\delta}^k| D^{\Phi_i(\bm{\lambda})}\langle r-2k\rangle$.
\item[2.]  $(i\text{-}\mathrm{Ind})^{(\delta)}\ D^{\Phi_i(\bm{\lambda})} \sim \sum\limits_{k=0}^{r} |\mathfrak{S}_{\delta}^k| D^{\bm{\lambda}}\langle r-2k\rangle.$
\item[3.] The blocks $B$ and $\Phi_i(B)$ have the same graded decomposition matrix.

\end{enumerate}
\end{thm}

\section{Categorical Actions and Crystals}
\subsection{Kashiwara Crystals}
We take a moment to recall the definition of crystals.

\begin{defn} (\cite[Definition 2.13]{BumpSchilling}) Given a root system $\Phi$ indexed by $I$ and a finite set $\mathcal{B}$, we can define a \textbf{crystal} via the following maps:
    \begin{align*}
        \tilde{e}_i, \tilde{f}_i &: \mathcal{B} \rightarrow \mathcal{B} \sqcup \{ 0\} \\
        \varepsilon, \varphi &: \mathcal{B} \rightarrow \mathbb{Z} \sqcup \{ - \infty \} \\
        \operatorname{wt} &: \mathcal{B} \rightarrow X
    \end{align*}
    such that the following holds:
   \[ \tilde{e}_i (x) = y\iff \tilde{f}_i(y) = x\]
    in which case:
    \begin{align*}
        \varepsilon_i(y) = \varepsilon_i(x)-1 \\
        \varphi_i(x)= \varphi_i(y)-1 \\
      \langle \operatorname{wt}(x), \alpha_i^\vee \rangle = \varphi_i(x) - \varepsilon_i(x)
    \end{align*} 
\end{defn}
We can visualize such a structure as a directed graph with vertices indexed by $\mathcal{B}$, and arrows labelled by $i \in I$. $\tilde{f}_i(x) = y$ then means we have an $i$-arrow from node $x$ to node $y$, whereas $\tilde{e}_i(y)=x$ can be seen as traversing the $i$-arrow $x \rightarrow y$ backwards. 
\[ x \overset{i}{\longrightarrow}y \quad \sim \quad \tilde{f}_i(x)=y, \quad \tilde{e}_i(y)=x\]
Given $i \in I$, not all nodes have an ingoing or outgoing $i$-arrow, and each node has at most one outgoing $i$-arrow and one incoming $i$-arrow (we say that $\tilde{f}_i(x)=0$ when there is no outgoing $i$-arrow from $x$, and $\tilde{e}_i(x) = 0$ if there is no incoming $i$-arrow). The constraints on the functions $\varepsilon_i$ and $\varphi_i$ mean that, given a node $x$ in the graph, $\varepsilon_i(x)$ gives the number of times it is possible to move ``backward" before there is not an incoming $i$-arrow to traverse (how many times $\tilde{e}_i$ can be applied before obtaining $0 \notin \mathcal{B}$) and $\varphi_i(x)$ counts the number of times it is possible to apply $\tilde{f}_i$ before obtaining $0$. More detailed exposition on this topic can be found in \cite{BumpSchilling}, but we want to emphasize that this is a fundamentally combinatorial object with a visual presentation. Since we will be working with Ariki--Koike algebras, we will mainly be considering crystals of type $\widehat{A}_{e-1}$.

A key component of crystals that we will require is their root strings.

\begin{defn}[\cite{BumpSchilling} \textsection 2.5]
    Given $x, y \in \mathcal{B}$ we say that $x \sim y$ if there exists some integer $k$ such that $x = \tilde{e}_i^k(y)$ or $x=\tilde{f}_i^k(y)$. This partitions the crystal into equivalence classes, where if $x \sim y$ we say that $x$ and $y$ belong to the same \textbf{root string}.
\end{defn}

We note that:
\begin{itemize}
    \item The weights of the elements of an $i$-root string differ by integers multiples of $\alpha_i$.
    \item Given a particular root string in a crystal, there is a unique element $b_h$ in this specific string such that $\tilde{e}_i(b_h)=0$ (and similarly, a unique element $b_l$ in the string such that $\tilde{f}_i(b_l)=0$). Visually, we can isolate the sub-graph given by this root string as:
    \[ \begin{tikzpicture}
        \filldraw (0,0) circle (2pt);
        \filldraw (1,0) circle (2pt); 
         \filldraw (5,0) circle (2pt);
         \node at (3,0) {\dots};
         \draw[thick, ->] (0,0) -- node[above] {$\tilde{f_i}$} (1,0);
          \draw[thick, ->] (1,0) -- node[above] {$\tilde{f_i}$} (2,0);
           \draw[thick, ->] (4,0) -- node[above] {$\tilde{f_i}$} (5,0);
           \node[anchor=east] at (0,0) {$b_h$};
            \node[anchor=west] at (5,0) {$b_l$};
    \end{tikzpicture}.\]
\end{itemize}

\begin{defn}[\cite{BumpSchilling} Definition 2.35, Proposition 2.36]\label{crystalrefl}
    Given a node $b$ in a crystal, if $k=\langle \operatorname{wt}(b), \alpha_i^\vee\rangle$, then
    \[ \sigma_i(b) = \begin{cases}
        \tilde{f}_i^k(b) &\text{if } k >0\\
        b & \text{if } k = 0\\
        \tilde{e}_i^k(b) & \text{if } k<0
    \end{cases}\] defines a map of order two from the crystal to itself which preserves root strings, and such that for any $b \in \mathcal{B}$,
    \[ \operatorname{wt}(\sigma_i(b)) = s_i(\operatorname{wt}(b))\]
    where $s_i(\mu) = \mu - \langle \mu, \alpha_i^\vee \rangle \alpha_i$ is the usual reflection of the weight lattice in the hyperplane orthogonal to the root $\alpha_i$.
\end{defn}

\subsection{Categorical Actions}

Let $\mathfrak{g}$ be a Kac--Moody algebra with weight lattice $P$ and Cartan subalgebra $\mathfrak{h}$. Let $\{ \alpha_i\}_{i \in I}\subset \mathfrak{h}^*$ denote the simple roots. Khovanov--Lauda and Rouquier (\cite{KLI}, \cite{R}) gave independent constructions (which were proven to be equivalent by Brundan \cite{Br15}) of a $2$-category $\mathcal{U}(\mathfrak{g})$ such that:
\begin{itemize}
    \item object set $P$
    \item  $1$-morphisms generated by $\{ \Ei, \Fi\}_{i \in I}$ acting on $P$ as:
    \begin{equation}\label{mathcalops}
       \Ei: \lambda \rightarrow \lambda + \alpha_i \qquad \Fi: \lambda \rightarrow \lambda - \alpha_i 
    \end{equation}
    \item $2$-morphisms as described in \cite{Br15} (we will not require their details here).
\end{itemize}

\begin{defn}
    A \textbf{categorical action} of $\mathfrak{g}$ is a representation of $\mathcal{U}(\mathfrak{g})$ defined above. In particular, we assign a category $\mathcal{C}_\nu$ to each weight $\nu \in P$, a functor to each of $\Ei$ and $\Fi$, and natural transformations to each $2$-morphism which descend to the NilHecke relations in the Grothendieck group upon decategorification.

\end{defn}

In other words, we have a decomposition of a category $\mathcal{C} = \bigoplus_{\nu \in P} \mathcal{C}_\nu$, with endofunctors $\Ei$ and $\Fi$ acting in the same way as the Chevalley generators.

Brundan and Davidson showed that this can be encoded combinatorially using crystals, where (for appropriate choices of types and parameters) the Kashiwara crystal operators $\tilde{e}_i$ and $\tilde{f}_i$ correspond to the functors $\Ei$ and $\Fi$, respectively.

\begin{thm}[\cite{BD} Theorem 4.31]\label{crystaliso}
    Given a nilpotent categorical action on a locally Schurian category with simple objects $\{ L(b) \mid b \in B\}$ \footnote{These definitions can be found in \cite{BD}, but the reader can rest assured that $\bigoplus_n\ak{\Lambda}{n}\mmod$ with simples $\simple{\ak{}{}}{\lambda}$ satisfies them}, there is a unique crystal structure on $B:= \sqcup_{\lambda} B_\lambda$ such that the following hold:
    \begin{enumerate}
        \item $\tilde{e}_ib \neq 0 $ if and only if $\Ei (L(b))\neq 0$; if this is the case, then \[\operatorname{soc}(\Ei(L(b)) \cong \operatorname{hd}(\Ei L(b)) \cong L(\tilde{e}_ib)\]
          \item $\tilde{f}_ib \neq 0 $ if and only if $\Fi (L(b))\neq 0$  if this is the case, then \[\operatorname{soc}(\Fi(L(b)) \cong \operatorname{hd}(\Fi L(b)) \cong L(\tilde{f}_ib)\]
          \item $\operatorname{soc}(\Ei^{(n)} (L(b)) \cong \operatorname{hd} (\Ei^{(n)}(L(b)) \cong L(\tilde{e}_i^nb).$
           \item $\operatorname{soc}(\Fi^{(n)} (L(b))) \cong \operatorname{hd} (\Fi^{(n)}(L(b)) \cong L(\tilde{f}_i^nb).$
    \end{enumerate}
    
\end{thm}

That is to say, the crystal encodes the structure of the simples in the module category. Each vertex $b \in \mathcal{B}$ indexes a unique simple module; applying a weight-lowering operator $\tilde{f}_i$ in the crystal gives the vertex indexing the unique simple socle (equivalently head) of the result of applying the functor $\Fi$ to the weight space indexed by the vertex $b$ (and analogously for weight-raising operators). To compensate for the occurrence of repeated summands when iterating these functors, the repeated crystal operator $\tilde{e}_i^n$ corresponds to the divided difference operator $\Ei^{(n)}$ on the side of the module category.

In \cite{CR}, Chuang and Rouquier initially developed the machinery of a categorical $\mathfrak{sl}_2$-representation and used it to establish Broué's abelian defect group conjecture: any two blocks of modules over $\mathbb{F}_p \mathfrak{S}_n$ and $\mathbb{F}_p \mathfrak{S}_m$ with the same defect group are derived equivalent. 
 
In other words, given $\mathcal{C} = \bigoplus_{n \in \mathbb{Z}} \mathcal{C}_n$ with $\mathcal{E}: \mathcal{C}_n \rightarrow \mathcal{C}_{n+2}$ and $\mathcal{F}: \mathcal{C}_n \rightarrow \mathcal{C}_{n-2}$ (with natural transformations inducing the $\mathfrak{sl}_2$ relations at the Grothendieck group-level), one obtains
\[ \mathcal{D}^b(\mathcal{C}_n) \cong \mathcal{D}^b(\mathcal{C}_{s_in})\]
via the categorification of the action of the reflection in the Weyl group. This equivalence is, in general, perverse: it acts on the derived category (the category of chain complexes on modules up to homotopy), sending simples to shifts of simples (possibly modulo lower-order terms); it does not send simples to simples. However, when the perversity function is $0$, the equivalence sends simples to simples and is $t$-exact.

In \cite{rock}, Webster interpreted the derived equivalences of \cite{CR} in the general setting (that is, $\mathfrak{g}$ acting on $\mathcal{C} = \bigoplus_{\nu \in P} \mathcal{C}_\nu$, as in \cref{mathcalops}):

\begin{thm}[\cite{rock} Lemma 3.1]\label{rockmain}
Let $\alpha_i^\vee( \nu) = k >0$. The divided difference functor $\Fi^{(k)}: \Cat{\nu} \rightarrow \Cat{s_i \nu}$ is a Morita equivalence if and only if $\Cat{\nu+\alpha_i}=0.$
\end{thm}

Therefore, in fact, \textit{most} Chuang--Rouquier equivalences are $t$-exact; Webster called these \textbf{Scopes equivalences} as they recover the Morita equivalences between blocks of $\mathbb{F}\mathfrak{S}_n$-modules established by Scopes in \cite{scopesCartanMatrices1991}.

Since they arise from $t$-exact Chuang--Rouquier equivalences, they are not only equivalences of abelian categories: they induce bijections between simples in each weight space. In the case of a cellularity-preserving idempotent truncation of a quasi-hereditary algebra, this further gives us equivalences of decomposition matrices --- this is the case of Ariki--Koike algebras, which have $q$-Schur algebras as constructed in \cite{djm} as quasi-hereditary covers.

\section{Categorical Actions on $\bigoplus_n\ak{\Lambda}{n}\mmod$}
We recall that we have defined $\ak{\Lambda}{n}$ with parameter $q$ in the quadratic relation such that $e$ was minimal satisfying $1+q+\dots +q^{e-1}=0$ (i.e. of quantum characteristic $e$).

\begin{thm}[\cite{BKgr} \S 4]\label{isommodule}
There is an isomorphism between an integrable highest weight module $V(\Lambda)$ of $U_q(\widehat{\mathfrak{sl}}_e)$ and the Grothendieck group of the category of modules over $\ak{\Lambda}{n}$, where the Chevalley generators $e_i$ and $f_i$ correspond to $i$-induction and $i$-restriction on graded modules of $\ak{\Lambda}{n}$.
\end{thm}

That is, we obtain a categorical $\widehat{\mathfrak{sl}}_e$-action on $\ak{\Lambda}{n}\mmod$, with the weight space decomposition given by the blocks of the Ariki--Koike algebra: we assign a block $B_\alpha$ to $\mathcal{C}_{\Lambda-\alpha}$, where $\alpha=\alpha(\bm\la)=\sum_i c_i(\bm\la) \alpha_i$ encodes the residue set of a multipartition $\bm\la$ in the block, see \cref{AK}.
Note that given a simple $D^{\bm\lambda} \in B_{\alpha}$, 
\begin{itemize}
    \item  $\tilde{e}_i$ removes a single $i$-node from the multipartition:
    $\tilde{e}_i([\bm\la])=[\bm\la]\setminus \{\mathfrak{n}\}$ where $\{\mathfrak{n}\}$ is a removable $i$-node of $\bm\la$.
    The corresponding weight space change is \[\mathcal{E}_i(\mathcal{C}_{\Lambda-\alpha} )  = \mathcal{C}_{\Lambda -( \alpha - \alpha_i)}= \mathcal{C}_{(\Lambda - \alpha) + \alpha_i}.\]
    That is, we \textit{add} $\alpha_i$ to the weight when applying $i$-restriction.
    \item similarly, $\tilde{f}_i$ adds a single $i$-node to the multipartition: $\tilde{f}_i([\bm \lambda]) = [\bm \lambda]\cup \{ \mathfrak{n}\}$ where $\mathfrak{n}$ is an addable $i$-node. The corresponding weight space change is \[\mathcal{F}_i(\mathcal{C}_{\Lambda-\alpha})  = \mathcal{C}_{\Lambda -( \alpha + \alpha_i)}= \mathcal{C}_{(\Lambda - \alpha) - \alpha_i}.\]
    That is, we \textit{subtract} $\alpha_i$ from the weight when applying $i$-induction.
\end{itemize}

Note that in the table below, in the case of the Ariki--Koike algebra $\ak{\Lambda}{n}$, $\nu = \Lambda - \alpha(\bm \lambda).$

\begin{figure}[h]
    \centering
    \scalebox{0.75}{
    \begin{tabular}{c|c|c|c|c}
    Crystal & Categorical Module & $\ak{\Lambda}{n}\mmod$ & Abacus & Multi-partition \\ \hline
    $\tilde{e}_i$ & $\mathcal{E}_i$ & $i-\operatorname{Res}$ & push a bead from runner $i$ to $i-1$ & remove an $i$-node\\
     $\tilde{f}_i$ & $\mathcal{F}_i$ & $i-\operatorname{Ind}$ & push a bead from runner $i-1$ to $i$ & add an $i$-node\\
    $\tilde{e}_i^k$ & $\mathcal{E}_i^{(k)}$ & $(i-\operatorname{Res})^{(k)}$
 & swap runners $i-1$ and $i$ & remove all removable $i$-nodes \\
  $\tilde{f}_i^k$ & $\mathcal{F}_i^{(k)}$ & $(i-\operatorname{Ind})^{(k)}$
 & swap runners $i-1$ and $i$ & add all addable $i$-nodes\\
 $\varepsilon_i(\nu)=0$ & $\mathcal{C}_{\nu + \alpha_i} = 0$ & $(i-\operatorname{Res}(D^\nu)) = 0$ & bead condition in \cref{AliceMain} & no removable $i$-nodes
 \\
 $\varphi_i(\nu)=0$ & $\mathcal{C}_{\nu - \alpha_i} = 0$ & $(i-\operatorname{Ind}(D^\nu)) = 0$ & bead condition in \cref{AliceMain} & no addable $i$-nodes
    \end{tabular}}

    \caption{Translating between Crystals, Modules, and Combinatorics}
    \label{fig:placeholder}
\end{figure}

The discussion in \cite[\textsection 3]{rock} applies directly to the case of the Ariki--Koike algebra with quasi--hereditary cover given by the cyclotomic $q$-Schur algebra:

\begin{defn}[\cite{djm} Definition 6.1]
 The cyclotomic $q$-Schur algebra is defined (up to Morita equivalence) as \[ \mathcal{S}(n) = \operatorname{End}_{\ak{\Lambda}{n}} \left( \bigoplus_{\mu} M^\mu\right)\] 
   where the $M^\mu$ are a collection of permutation modules indexed by $\ell$-multicompositions (see \cite[\textsection 3]{djm} for their construction). Let $\mathcal{S} := \bigoplus_n \mathcal{S}(n).$
\end{defn}

\begin{thm}[\cite{djm} Theorem 6.16, Corollary 6.18]\label{cyclotomicschur}\;
The cyclotomic $q$-Schur algebra is cellular, with standard modules $\std{\mathcal{S}}{\bm \mu}$ (equivalently, simples $\simple{\mathcal{S}}{\bm \mu}$) indexed by $\ell$-multipartitions $\bm\mu = (\mu^{(1)}, \dotsc, \mu^{(\ell)}).$
    
\end{thm}

By construction, $\mathcal{S}\mmod$ is highest weight, and the Schur functor induces a quotient to the category of Ariki--Koike modules.

\begin{prop}[\cite{rock} Proposition 3.4]
Let $\mathcal{C} =\mathcal{S}\mmod$. Then $\overline{\mathcal{C}}=\bigoplus_n\ak{\Lambda}{n}\mmod$ can be obtained as a quotient of $\mathcal{C}$ via the Schur functor: we mod out by the Serre subcategory $\mathcal{B} \subset \mathcal{C}$ of $\mathcal{S}$-modules $M$ such that $M\bm e = 0$ (where $\bm e$ is the idempotent such that $\ak{\Lambda}{n} = \bm e \mathcal{S}(n) \bm e$). If $\nu= \Lambda - \alpha(\bm\lambda) \in P$ and $i \in I$ are chosen so that \cref{rockmain} holds, then:

\begin{itemize}
    \item $\Fi^{(k)}: \mathcal{C}_{\nu} \rightarrow \mathcal{C}_{s_i \nu}$ is an equivalence of highest weight categories sending standard modules to standard modules $\std{\mathcal{S}}{\bm \lambda} \mapsto \std{\mathcal{S}}{\tilde{f}_i^k(\bm \lambda)} \in \mathcal{S}\mmod$ and Specht modules to Specht modules: $S^{\bm \lambda} \mapsto S^{\tilde{f}_i^k(\bm \lambda)} \in \ak{\Lambda}{n}\mmod$,
    \item $\Fi^{(k)}$ maintains decomposition numbers of both the cyclotomic $q$-Schur algebra and the Ariki--Koike algebra,
\end{itemize}
where $k = \langle \nu, \alpha_i^\vee \rangle >0$ (and similarly using $\mathcal{E}_i^{(k)}$, $\tilde{e}_i^k$ when $k<0$).
\end{prop}

In terms of the crystal graph, \cref{rockmain} translates to having a Morita equivalence given by reflecting a root string precisely when the simple is at the very end of a root string --- either $\varphi_i(\nu) = 0$ or $\varepsilon_i(\nu) = 0$. In fact, this is the situation the Chuang--Rouquier construction categorifies the reflection of the root string \cref{crystalrefl} given in \cite[\textsection 2.5]{BumpSchilling}.

Looking at the combinatorial interpretation of the crystal rule in terms of (multi-)partitions, this occurs when we have no addable (or no removable) $i$-nodes.

In \cite{DA24}, the first author works with the established crystal action given by moving beads on abacus runners, as described in \cref{betaAbacus}. The key to note is that the condition in \cref{AliceMain} corresponds precisely to this root-string requirement: swapping two runners is equivalent to adding all addable and removing all removable $i$-nodes. If there were only nodes of one type (addable or removable), then the swap would only add (or only remove) nodes, as the other option would be impossible. As a result, trying to apply a single crystal operator of this second type (for example, trying to add a node when there are no addable nodes) would output a zero. By \cref{crystaliso} this is exactly Webster's condition for a $t$-exact Chuang--Rouquier equivalence. Knowing that these particular Morita equivalences preserve decomposition numbers (see \cite[Proposition 3.4]{rock}) is in agreement with the first author's study of the images of Specht modules under the associated maps.

\begin{cor}
The equalities between decomposition matrices of Ariki--Koike algebras established by the first author in \cref{AliceMain} are induced by the cases of Scopes equivalences when treating $\ak{\Lambda}{n}\mmod$ as a categorical $\widehat{\mathfrak{sl}}_e$-module.
That is:

\[ \begin{tikzpicture}
\node at (-1.5,0) {$\Fi^{(k)}:$};
\node at (-0.75,-1.5) {$\Phi_i:$};
    \node at (0.5,0) {$ \mathcal{C}_{\Lambda - \alpha}$};
    \node at (5,0) {$\mathcal{C}_{\Lambda - \alpha-k\alpha_i}$};
    \node at (0.5,-1.5) {$B_{\alpha}$};
    \node at (5,-1.5) {$B_{\Phi_i(\alpha)}$};
    \draw[->] (2,0)-- (3.5,0);
     \draw[->] (2,-1.5)-- (3.5,-1.5);
    \node at (0.5,-0.75) {$\parallel$};
     \node at (5,-0.75) {$\parallel$};

\end{tikzpicture}\]

where
\begin{itemize}
    \item $\alpha=\alpha(\bm\la)=\sum_i c_i(\bm\la) \alpha_i$ for $\bm\la$ a multipartition in the block $B_\alpha$;
    \item $\Phi_i(\alpha)$ encodes the residue set of any multipartition in the block of $\Phi_i(\bm \lambda)$.
    
\end{itemize} 
\end{cor}

\printbibliography

\end{document}